\documentclass{amsart}
\usepackage{amssymb, amsmath, amsfonts, amscd}

\input xypic

\theoremstyle{plain}
\newtheorem{theorem}{Theorem}[section]
\newtheorem{corollary}[theorem]{Corollary}
\newtheorem{lemma}[theorem]{Lemma}
\newtheorem{proposition}[theorem]{Proposition}

\newtheorem{example}[theorem]{Example}

\theoremstyle{definition}
\newtheorem{definition}[theorem]{Definition}

\theoremstyle{remark}

\numberwithin{equation}{theorem}

\newcommand{\A}{\mathcal{A}}

\newcommand{\E}{\mathcal{E}}
\renewcommand{\O}{\mathcal{O} }

\renewcommand{\P}{\mathbf{P} }
\newcommand{\Z}{\mathbf{Z} }
\renewcommand{\Pr}{\mathcal{P} }
\renewcommand{\u}{\frac{1}{u} }
\renewcommand{\t}{\frac{1}{t} }

\newcommand{\U}{\operatorname{U}}

\newcommand{\Proj}{\operatorname{Proj} }
\newcommand{\D}{\operatorname{D} }

\begin{document}

\title{Modules of principal parts on the projective line }
\author{Helge Maakestad }
\address{Department of Mathematics, Faculty of Exact Sciences, Bar-Ilan University,
 Ramat-Gan, Israel }
\email{makesth@macs.biu.ac.il }
\thanks{Research supported by the EAGER foundation, the Emmy Noether Research 
Institute for Mathematics, the Minerva Foundation of Germany and the 
Excellency Center ``Group Theoretic Methods in the Study of Algebraic 
Varieties''}
\date{2001}
\begin{abstract} 
The modules of principal parts $\Pr^k(\E)$ of a locally free sheaf $\E$
on a smooth scheme $X$ is a sheaf of $\O_X$-bimodules which is locally free as left and right 
$\O_X$-module. We split explicitly the modules of principal parts
$\Pr^k(\O(n))$ on the projective line in arbitrary characteristic, as left
and right $\O_{\P^1}$-module. We get examples when the splitting-type
as left module differs from the splitting-type as right module. We also
give examples showing that the splitting-type of the principal parts
changes with the characteristic of the base field.
\end{abstract}
\maketitle

\tableofcontents

\section{Introduction}
In this paper we will study the splitting-type of the modules of 
principal parts
of invertible sheaves on the projective line as left and right 
$\O_{\P^1}$-module in 
arbitrary characteristic. 
The splitting type 
of the principal parts $\Pr^k(\O(n))$ as a left $\O_{\P^1}$-module 
in  characteristic zero has been studied by several authors 
(see \cite{PERK2}, 
\cite{PIE} and \cite{dirocco}).
The novelty of this work is that  we consider the principal parts as left 
and right $\O_{\P^1}$-module in arbitrary characteristic.  We  
give examples when the splitting-type as left $\O_{\P^1}$-module differs 
from the splitting-type
as right $\O_{\P^1}$-module. 
The main theorem of the paper (Theorem \ref{leftandright}),
gives the splitting-type of $\Pr^1(\O(n))$  as left and right 
$\O_{\P^1}$-module
for all $n\geq 1$ over any field $F$. The result is the following: 
The principal  parts $\Pr ^1(\O (n))$ splits as $\O (n)\oplus \O (n-2)$ 
as right $\O _{\P^1}$-module. If the characteristic of 
$F$ divides $n$, $\Pr ^1(\O (n))$ splits as $\O(n)\oplus \O(n-2)$ as left 
$\O_{\P^1}$-module. If the characteristic of $F$ does not divide $n$, 
$\Pr^1(\O(n))$ splits as $\O(n-1)\oplus \O(n-1)$ as left $\O_{\P^1}$-module.
Hence the modules of principal parts are the first 
examples of a sheaf of abelian groups equipped with two non-isomorphic 
structures  as locally free sheaf. 
In the papers \cite{PERK2}, \cite{PIE} and \cite{dirocco} the authors work
in characteristic zero, and they only consider the left
module structure. The proofs are not explicit. In this work we split
the principal parts explicitly and the techniques we develop will be used in 
future papers to get deeper knowledge of the principal parts in positive
characteristic. 
All computations are explicit, and in sections 3 - 6 we 
develop techniques to construct explicit non-trivial maps of $\O$-modules
from $\O(l)$ to $\Pr^k(\O(n))$. The main theorem here is Theorem \ref{maintheorem},
where we prove existence of certain systems of linear equations with integer
coefficients. Solutions to the systems satisfying extra criteria determins
the splitting-type of $\Pr^k(\O(n))$. Theorem \ref{maintheorem}
is used in Proposition \ref{characteristiczero} to determine the splitting-type
of $\Pr^k(\O(n))$ for all $1\leq k \leq n$ in characteristic zero, and we recover
results obtained in  \cite{PERK2},\cite{PIE} and \cite{dirocco}.
We also give examples where the splitting-type
can be determined by diagonalizing the structure-matrix defining the 
principal parts (Section 4).

\textbf{Acknowledgments} This paper form part of my PhD-thesis written 
under supervision of Dan Laksov at the Royal Institute of Technology 
in Stockholm, and I want to thank him for expert advice on algebraic 
geometry.  
The paper in its final form was 
written partly in spring 2001 during a stay at the 
Institut Joseph Fourier 
in Grenoble, financed by an EAGER-scholarship, and I want to thank Michel 
Brion for an invitation to the IJF. Parts of the paper was written in 
September 2001 during a stay at the Department of Mathematics at the Bar-Ilan 
University in Israel 
and I want to thank Mina Teicher for an invitation to the Bar-Ilan University.
Thanks also to  Rolf K\"{a}llstr\"{o}m
for interest in this work.

\section{Modules of principal parts}
In this section we define and prove basic properties of the principal parts:
existence of fundamental exact sequences, functoriality and existence
of bimodule-structure. In the following let $X$ be a scheme defined over 
a fixed
base-scheme $S$. We assume that $X$ is separated and smooth over $S$. Let $\Delta$
in $X\times _S X$ be the diagonal, and let $\mathcal{I}$ in $\O_{X\times_S X}$ be the sheaf of ideals defining $\Delta$. Let $X^k$ be the scheme with topological space
$\Delta$ and with structure-sheaf $\O_{\Delta^k}=\O_{X\times _S X}/\mathcal{I}^{k+1}$.
By definition, $X^k$ is the \emph{$k$'th order infinitesimal neighborhood of the diagonal}. In the following, we omit reference to the base scheme $S$ in products.
\begin{definition} \label{principalparts} 
Let $\E$ be an quasi-coherent $\O_X$-module. We define
the \emph{$k$'th order modules of principal parts} of $\E$, to be
\[ \Pr^k_X(\E)=p_*(\O_{\Delta^k}\otimes q^*\E).\]
We write $\Pr^k_X$ for the module $\Pr^k_X(\O_X)$.
\end{definition}
When it is clear from the context which scheme we are working on, we 
write $\Pr^k(\E)$ instead of $\Pr^k_X(\E)$. 

\begin{proposition} \label{fundamentalexact} Let $\E$ be a quasi-coherent
$\O_X$-module.
There exists an exact sequence
\[ 0 \rightarrow S^k(\Omega^1_X)\otimes \E \rightarrow \Pr^k_X(\E) \rightarrow
\Pr^{k-1}_X(\E) \rightarrow 0 \]
of $\O_X$-modules, where $k=1,2,\dots$.
\end{proposition}
\begin{proof} See \cite{LAKSOV}, Section  4.\end{proof}
From Proposition \ref{fundamentalexact} 
it follows by induction, that for a locally free sheaf $\E$ of rank 
$e$, $\Pr^k_X(\E)$ is locally free of rank $e\binom{n+k}{n}$, where $n$ is the 
relative dimension of $X$ over $S$. 
\begin{proposition} \label{functoriality}
Let $f:X\rightarrow Y$ be a map of smooth schemes
over $S$, and let 
$\E$ be a locally free $\O_Y$-module. 
Then there exists a commutative diagram of exact sequences
\[ \diagram 0 \rto & S^k(f^*\Omega^1_Y)\otimes f^*\E \dto \rto
& f^*\Pr^k_Y(\E) \dto \rto & f^*\Pr^{k-1}_Y(\E) \dto \rto & 0 \\
0 \rto & S^k(\Omega^1_X)\otimes f^*\E \rto & \Pr^k_X(f^*\E) \rto & \Pr^{k-1}_X(f^*\E)
\rto & 0 \enddiagram \]
of left $\O_X$-modules for all $k=1,2,\dots$.
\end{proposition}
\begin{proof} See \cite{PERK}.
\end{proof}
From Proposition \ref{functoriality} it follows that if $U$ is an open subset of $X$, we have that
$\Pr^k_X(\E)|_U$ equals $\Pr^k_U(\E|_U)$, which allows one to do local 
computations.
\begin{proposition} \label{functor}
The principal parts $\Pr^k_X$ define
a covariant functor
\[ \Pr^k_X: \O_X-\underline{mod}\rightarrow \Pr^k_X-\underline{mod}, \]
where for all quasi-coherent $\O_X$-modules $\E$, the k'th order principal parts
$\Pr^k_X(\E)$ is a quasi-coherent $\Pr^k_X$-module. The functor is right exact and commute with direct limits. If $\Pr^k_X$ is flat, the functor is exact.
\end{proposition}
\begin{proof} See \cite{EGAIV}, Proposition 16.7.3.\end{proof}

Note that since we assume $X$ smooth over $S$, it follows that 
$\Pr^k_X$  is locally free, hence the functor in Proposition 
\ref{functor} is exact.

We next consider the bimodule-structure of the principal parts.
%
\begin{proposition}\label{pushiso} 
Let $f,g:U\rightarrow V$ be morphisms of topological spaces,
and let $s$ be a section of $f$ and $g$ with $s(V)$ a closed set. Let 
furthermore $\A$ be a sheaf of abelian groups on $U$ with support in $s(V)$
Then $f_*(\A)$ equals $g_*(\A)$.
\end{proposition}
\begin{proof} We first claim that the natural map from $s_*s^{-1}\A$
to $\A$ is an isomorphism: Since $s(V)$ is closed, and $s$ is a section of $f$, 
it follows that $s$ is a closed map. Both $s_*s^{-1}\A$ and $\A$ have support 
contained in $s(V)$, hence we prove that the map is an isomorphism at the stalks for 
all points $p$ in $s(V)$. We see that $(s_*s^{-1}\A)_{s(p)}$ equals
$(s^{-1}\A)_p$, since $s$ is a closed immersion. Furthermore we have that 
$(s^{-1}\A)_p$ equals $\A_{s(p)}$ since $s^{-1}$ is an exact functor, 
and the claim follows. It now follows that $f_*\A$ is isomorphic to
$f_*s_*s^{-1}\A$, and since $s$ is a section of $f$, we get that
$f_*\A$ is isomorphic to $s^{-1}\A$. A similar argument proves that
$g_*\A$ is isomorphic to $s^{-1}\A$, and the proposition follows.
\end{proof}

Let now $\Delta$ be the diagonal in $X\times X$, which is closed since
$X$ is separated over $S$. Then the sheaf 
$\O_{\Delta^k}\otimes q^*\E$ has support in $\Delta $, hence by Proposition
\ref{pushiso} it follows that $\Pr^k_X(\E)=p_*(\O_{\Delta^k}\otimes q^*\E)$ is 
isomorphic to $q_*(\O_{\Delta^k}\otimes q^*\E)$. By the projection formula,
we see that $q_*(\O_{\Delta^k}\otimes q^*\E)$ equals 
$q_*(\O_{\Delta^k})\otimes \E$, hence by Proposition \ref{pushiso}
we see that $\Pr^k_X(\E)$ is isomorphic to
$q_*(\O_{\Delta^k})\otimes \E$
as sheaves of abelian groups. Identifying $q_*(\O_{\Delta^k})\otimes \E$ with
$\Pr^k_X(\E)$, we see that we have defined two $\O_X$-module structures on 
$\Pr^k_X(\E)$, hence it is a sheaf of $\O_X$-bimodules. This means that
for any open set $U$ of $X$, the sheaf $\Pr^k_U(\E|_U)$ 
is an $\O_X(U)$-bimodule, and all restriction maps are maps of bimodules
satisfying obvious compatibility criteria.
Let $X^k$ be the $k$'th order infinitesimal neighborhood of the diagonal.
Then the two projection-maps $p,q:X^k\rightarrow X$ 
induce two maps $l,r:\O_X\rightarrow \Pr^k_X$ of $\O_X$-modules. The maps
$l$ and $r$ are the maps defining the bimodule structure on $\Pr^k_X$, and
we see that $\Pr^k_X$ is a sheaf of $\O_X$-bialgebras. The map
$d=l-r:\O_X\rightarrow \Pr^k_X$ is verified to be a differential operator
of order $k$, called \emph{the universal differential operator}. 

\section{Transition matrices for principal parts as left module}
In this section we compute explicitly the transition-matrices defining the principal
parts $\Pr^k(\O(n))$ on the projective line over the integers.
%
%
%
We will use the following notation: Define $\P^1$ as $\Proj\mathbf{Z}[x_0,x_1]$, 
where $\mathbf{Z}$ are the integers and put $\U _i=\D (x_i)$ and $\U _{01}=\D (x_0x_1)$ 
for $i=0,1$, where $x_i$ are homogeneous coordinates on $\P ^1$. 
Consider the modules of principal parts $\Pr^k$ from Definition 
\ref{principalparts} on $\P^1$ for $k\geq 1$.
On the open set $\U_{01}$, the modules of principal parts $\Pr ^k$ equals   
$\Z[t,\t,u,\u]/(u-t,\u-\t)^{k+1}$ as $\O_{\U_{01}}$-module, and $\O_{\U_{01}}$
is isomorphic to $\Z[t,\t]$. 

\begin{lemma} On the open set $\U_{01}$, as a left $\O$-module, 
$\Pr^k$ is a free  $\Z[t,\t]$-module of rank $k+1$, and there exists 
two natural bases. The bases are
$B=\{1,dt,dt^2, \dots , dt^k\}$ and $B'=\{ 1,ds,ds^2, \dots , ds^k \}$, where
$dt^i=(u-t)^i$ and $ds^i=(\u-\t)^i$. 
\end{lemma}
\begin{proof} Easy calculation.
\end{proof}
\begin{proposition} Consider $\Pr^k$ as left $\O$-module on $\P^1$ 
with $k=1,2,\dots$. On the open set $\U_{01}$ the transition matrix 
$[L]^{B'}_{B}$ between
the two bases $B'$ and $B$ is given by the following formula
\[ ds^p= \sum_{i=0}^{k-p}(-1)^{i+p}\frac{1}{t^{i+2p} }
\binom{i+p-1}{p-1}dt^{i+p} \] for all $0\leq p \leq k$.
\end{proposition}
%
%
\begin{proof} By definition $ds^p$ equals $(\u-\t)^p$ in the 
module $\Z[t,\t,u,\u](u-t,\u-\t)^{k+1}$. It follows  that
\[ ds^p=(\u-\t)^p=\frac{(-1)^p(u-t)^p}{u^pt^p}=
(-1)^p\frac{1}{t^p}dt^p\frac{1}{u^p} ,\]
and since $u=t+u-t=t+dt$ we get 
\[ ds^p=(-1)^p\frac{1}{t^p}dt^p\frac{1}{(t+dt)^p}.\]
We see that 
\[ ds^p=(-1)^p\frac{1}{t^{2p} }dt^p\frac{1}{(1+dt/t)^p} ,\]
and using the identity
\[ \frac{1}{(1+\omega)^p}=\sum_{i\geq0}(-1)^i\binom{i+p-1}{p-1}\omega^i \]
we get 
\[ ds^p= (-1)^p\frac{1}{t^{2p }}dt^p\sum_{i\geq 0}(-1)^i\binom{i+p-1}{p-1}
(\frac{dt}{t})^i.\]
We put $dt^{k+1}=dt^{k+2}=\cdots =0$ and get
\[ ds^p=\sum_{i=0}^{k-p}(-1)^{i+p}\frac{1}{t^{i+2p}}\binom{i+p-1}{p-1}dt^{i+p}\]
and the proposition is proved.
\end{proof}

Consider the invertible $\O$-module $\O(n)$ on $\P ^1$, with $n\geq 1$. 
We want to study the principal parts $\Pr ^k(\O(n))$ with $1\leq k \leq n$
as a left $\O$-module.

\begin{lemma} \label{basis} 
On the open set $\U_{01}$ as left $\O$ module, $\Pr^k(\O(n))$
is a free $\Z[t,\t]$-module of rank $k+1$, 
and there exists two natural bases. The bases are
$C=\{1\otimes x_0^n,dt\otimes x_0^n,\dots , dt^k\otimes x_0^n\}$ and 
$C'=\{ 1\otimes x_1^n,ds\otimes x_1^n, \dots , ds^k\otimes x_1^n \}$, where
$dt^i=(u-t)^i$ and $ds^i=(\u-\t)^i$. 
\end{lemma}
\begin{proof} Easy calculation. \end{proof}

\begin{theorem} \label{transitionmatrix} Consider $\Pr^k(\O(n))$ as left
$\O$-module on $\P^1$. 
On the open set $\U_{01}$ the transition-matrix $[L]^{C'}_{C}$ 
between the bases $C$ and $C'$ is given by the following formula
\[ ds^p\otimes x_1^n= \sum_{i=0}^{k-p}(-1)^p\frac{1}{t^{i+2p-n}}\binom{n-p}{i}
dt^{i+p}\otimes x_0^n .\]
where $0\leq p \leq k$.
\end{theorem}
\begin{proof} By definition
\[ ds^p\otimes x_1^n= (\u-\t)^p\otimes t^nx_0^n= (-1)^p\frac{(u-t)^p}{u^pt^p}
u^n\otimes x_0^n.\]
Since $n-p\geq 0$ and $u=t+dt$ we get
\[ ds^p\otimes x_1^n=(-1)^p\frac{1}{t^p}dt^p(t+dt)^{n-p}\otimes x_0^n.\]
Using the binomial theorem, we get 
\[ ds^p\otimes x_1^n=(-1)^p\frac{1}{t^p}dt^p\sum_{i=0}^{n-p}\binom{n-p}{i}
t^{n-p-i}dt^i\otimes x_0^n .\]
By assumption $dt^{k+1}=dt^{k+2}=\cdots =0$, which gives
\[ds^p\otimes x_1^n= \sum_{i=0}^{k-p}(-1)^p\frac{1}{t^{i+2p-n}}\binom{n-p}{i}
dt^{i+p}\otimes x_0^n, \]
and the theorem follows.
\end{proof}

\begin{corollary} \label{locallyfree}$\Pr^k(\O (n))$ is locally free of rank $k+1$ on $\P ^1$ for all
$1\leq k \leq n$.
\end{corollary}
\begin{proof} We compute the determinant of the structure matrix of
$\Pr^k(\O (n))$, $[L]^{C'}_C$ from Theorem \ref{transitionmatrix}:
The matrix $[L]^{C'}_C$ is lower triangular, and it follows that the determinant
$|[L]^{C'}_C |$ equals
\[ \prod_{p=0}^k(-1)^p\frac{1}{t^{2p-n}}\binom{n-p}{0}. \]
We simplify to get
\[ |[L]^{C'}_C|=\prod_{p=0}^k(-1)^p\frac{1}{t^{2p-n} }=
(-1)^{kp}t^{\sum_{p=0}^k2p-n }=(-1)^{kp}t^{(n-k)(k+1)} \]
and we see that $|[L]^{C'}_C|$ is a unit in $\Z[t,\t]$. It follows that 
$[L]^{C'}_C$ is an invertible matrix, and the result follows.
\end{proof}

By Proposition \ref{functoriality}, formation of principal parts 
commute with base-extension, and it follows that
$\Pr^k(\O (n))$ is locally free of rank $k+1$ on $\P^1_A$ for any commutative
ring $A\neq 0$ with unit.

\begin{example} \label{transexample}By Theorem 
\ref{transitionmatrix} the transition matrix $[L]^{C'}_C$ for 
$\Pr^1(\O (n))$ is as follows:
\[ [L]^{C'}_C= 
 \begin{pmatrix} t^n & 0 \\
            nt^{n-1}  & -t^{n-2} 
 \end{pmatrix}
.\]
We compute the determinant $|[L]^{C'}_C|$ and find that it equals $-t^{2n-2}$.
\end{example}
\section{Splitting principal parts as left module by 
matrix-diagonalization}
In this section we will split explicitly the modules of principal parts.
We will work over $\P^1$ defined over $F$, where $F$ is a field.
By  \cite{GRO}, Theorem 2.1, we know that all locally free sheaves of finite rank on 
$\P ^1$ split into a direct 
sum of invertible sheaves, and we want to compute  explicitly the 
splitting-type for the sheaf
$\Pr ^1(\O (n))$ as left $\O$-module. From Lemma \ref{basis}
it follows that on the basic open set
$\U_0$, $\Pr^1(\O(n))$ is a free $F[t]$-module on the basis $C=\{1\otimes x_0^n,dt\otimes x_0^n\}$. On the open set $\U_1$, $\Pr^1(\O(n))$ is a free $F[s]$-module 
on the basis
$C'=\{1\otimes x_1^n, ds\otimes x_1^n\}$, where $s=\t$.
When we pass to the open set $\U_{01}=\U_0\cap \U_1$ we see that $\Pr^1(\O (n))$
has $C$ and $C'$ as bases as $F[t,s]$-module. On $\U_0$ consider the new
basis $D=\{1\otimes x_0^n,t\otimes x_0^n+ndt\otimes x_0^h\}$. Consider also the
new basis $D'=\{\t\otimes x_1^n+nds\otimes x_1^n, 1\otimes x_1^n \}$ on the open set 
$\U_1$. Notice that 
$D$ and $D'$ are bases if and only if the characteristic of $F$ does not 
divide $n$, hence let us assume this for the rest of this section. 
We first compute the base change matrix for $\Pr^1(\O (n))|_{\U_0}$ from $C$ to $D$,
and we get the matrix
\[ \text{[I]$^C_D$ }=
 \begin{pmatrix}   1   &   -\frac{1}{n}t \\
                   0   &   \frac{1}{n}
 \end{pmatrix}.
\]
We secondly compute the base change matrix for $\Pr^1(\O (n))|_{\U_1}$ from $D'$ 
to $C'$, and we get the matrix
\[ \text{ [I]$^{D'}_{C'}$ =}
  \begin{pmatrix}   \t     &     1  \\
                     n     &     0  
  \end{pmatrix}.
\]
In Example \ref{transexample} we saw that the transition matrix defining $\Pr^1(\O (n))$
is given by
\[ \text{[L]$^{C'}_C$=} 
 \begin{pmatrix} t^n & 0 \\
            nt^{n-1}  & -t^{n-2} 
 \end{pmatrix}.
\]
We see that if we let $D$ be a new basis for $\Pr^1(\O (n))$ as $F[t]$-module
on $\U_0$, and let $D'$ be a new basis for $\Pr^1(\O (n))$ as $F[s]$-module 
on $\U_1$, then the transition matrix $[L]^{D'}_{D}$ becomes
\[ \text{ [L]$^{D'}_D$=[I]$^C_D$[L]$^{C'}_C$[I]$^{D'}_{C'}$ } \]
which equals
\[ \begin{pmatrix}   1   &   -\frac{1}{n}t \\
                   0   &   \frac{1}{n}
 \end{pmatrix}
\begin{pmatrix}   t^n   &               \\
                   nt^{n-1}    &   -t^{n-2}
 \end{pmatrix}
 \begin{pmatrix}   \t   &    1 \\
                   n   &     0 
 \end{pmatrix}
.\]
It follows that
\[ \text{[L]}^{D'}_{D}=
\begin{pmatrix}   t^{n-1}  &   0 \\
                      0     &  t^{n-1} 
 \end{pmatrix}
,\]
hence as a left $\O$-module, the principal parts  $\Pr^1(\O (n))$ splits
as $\O (n-1)\oplus \O (n-1)$. By Proposition \ref{functoriality}, 
it follows that the splitting
$\Pr^1(\O (n))\cong \O(n-1)\oplus \O(n-1)$ is valid on $\P^1_A$, where $A$ is any
$F$-algebra.

\section{Maps of modules and systems of linear equations}
In this section we study criteria for splitting the principal parts
on the projective line $\P^1$ over any field $F$. 
Given $\Pr^k(\O(n))$ with $1\leq k \leq n$, we prove existence of 
systems of linear equations $\{A_r\mathbf{x_r}=b_r\}_{r=0}^k$, where
$A_r$ is  a rank $r+1$ matrix with integral coefficients. 
A solution $\mathbf{x_r}$ to the system $A_r\mathbf{x_r}=b_r$ 
gives rise to a map $\phi_r$  of left $\O$-modules from $\O(l)$ 
to $\Pr^k(\O(n))$. The main result is Theorem \ref{maintheorem}
where we prove that if there exists for all $r=0,\dots ,k$  
solutions $\mathbf{x_r}$ 
with coefficients in a field $F$, satisfying certain explicit 
criteria, then we can completely determine the splitting-type of the 
principal parts.

By  Corollary \ref{locallyfree} we know that $\Pr^k(\O (n))$ is 
locally free of rank $k+1$ 
over $\P^1$ defined over $\Z$, where $\Z$ is the integers, 
hence by base extension, $\Pr^k(\O (n))$ is locally free over $\P^1$ defined 
over any field $F$.
By \cite{GRO}, Theorem 2.1, we know that on $\P^1$ every locally 
free sheaf 
of finite rank splits uniquely into a direct sum of invertible 
$\O$-modules. 
Recall from Lemma \ref{basis} that $\Pr^k(\O (n))$ has two natural bases 
on the open set $\U_{01}$:
\[ C=\{1\otimes x_0^n,\dots , dt^k\otimes x_0^n \} \]
and
\[ C'=\{1\otimes x_1^n,\dots , ds^k\otimes x_1^n \} .\]
By theorem \ref{transitionmatrix}, the transition matrix $[L]^{C'}_C$ 
from basis $C'$ to $C$
is given by the following relation
\begin{equation} \label{basisrelation} ds^p\otimes x_1^n= 
\sum_{i=0}^{k-p}(-1)^p\frac{1}{t^{i+2p-n}}\binom{n-p}{i}
dt^{i+p}\otimes x_0^n .\end{equation}
We will use relation \ref{basisrelation} to construct split injective maps
$\O(l)\rightarrow \Pr^k(\O (n))$ of left $\O$-modules. Put $l=n-k$. On the open set
$\U_0$, the sheaf $\O(n-k)$ is isomorphic to $F[t]x_0^{n-k}$ as $\O$-module.
On the open set $\U_1$, $\O(n-k)$ is isomorphic to $F[\t]x_1^{n-k}$. For $i=0,1$ 
we want to  define maps 
\[\phi_r^i: \O (n-k)|_{\U_i}\rightarrow \Pr^k(\O (n))|_{\U_i} \]
where $r=0,\dots, k$, of left $\O_{\U_i}$-modules, agreeing on the open set 
$\U_{01}$. The maps 
$\{\phi_r^i\}_{i=0,1}$ will then glue to give $k+1$ well-defined maps 
$\phi_r:\O(n-k)\rightarrow \Pr^k(\O (n))$ of left $\O$-modules.
Define the map $\phi_0^1$ as follows: $\phi_0^1(x_1^{n-k})=1\otimes x_1^n$.
On the open set $\U_{01}$ we have that $x_0^{n-k}=t^{k-n}x_1^{n-k}$. We want to 
define $\phi^0_0(x_0^{n-k})$. Since $x_0^{n-k}=t^{k-n}x_1^{n-k}$ it 
follows that
$\phi^0_0(x_0^{n-k})=t^{k-n}\phi_0^1(x_1^{n-k})=t^{n-k}(1\otimes x_1^{n-k})$.
We now use relation \ref{basisrelation} which shows that
\[ 1\otimes x_1^n= \binom{n}{0}t^n\otimes x_0^n+\binom{n}{1}t^{n-1}dt\otimes x_0^n
+\cdots + \binom{n}{k}t^{n-k}dt^{n-k}\otimes x_0^n .\]
We define $\phi^0_0(x_0^{n-k})=t^{k-n}(1\otimes x_0^n)$, and it follows that
\[ \phi^0_0(x_0^{n-k})=t^{k-n}(\binom{n}{0}t^n\otimes x_0^n+\binom{n}{1}t^{n-1}dt\otimes x_0^n +\cdots + \binom{n}{k}t^{n-k}dt^{n-k}\otimes x_0^n) \]
which equals
\[ t^{k-n}\sum_{i=0}^k\binom{n}{i}t^{n-i}dt^i\otimes x_0^n.\]
Define $c^0_i=\binom{n}{i}$ and $x_{0,0}=1$. We get
\[ \phi^0_0(x_0^{n-k})=\sum_{i=0}^kc^0_it^{k-i}dt^i\otimes x_0^n \]
and 
\[ \phi_0^1(x_1^{n-k})=x_{0,0}(1\otimes x_1^n).\]
We see that we have defined a map of left $\O$-modules $\phi_0:\O(n-k)\rightarrow
\Pr^k(\O (n))$, which in fact is defined over the integers $\Z$.
We want to generalize the construction made above, and define maps of
left $\O$-modules $\phi_r:\O(n-k)\rightarrow \Pr^k(\O (n))$ for $r=1,\dots ,k$.
Define 
%
%
\begin{equation} \label{phi1def} \phi_r^1(x_1^{n-k})=x_{0,r}\frac{1}{t^r}
\otimes x_1^n+x_{1,r}\frac{1}{t^{r-1}}ds\otimes x_1^n+ \end{equation} 
\[
\cdots +x_{r-1,r}\t ds^{r-1}\otimes x_1^n+x_{r,r}ds^r\otimes x_1^n . \]
where the symbols $x_{i,r}$ are independent variables over $F$ for all $i$ and $r$.
Simplifying, we get
\[ \phi_r^1(x_1^{n-k})=\sum_{j=0}^rx_{j,r}t^{j-r}ds^j\otimes x_1^n. \]
We want to define a map $\phi_r^0$ on $\U_0$, such that $\phi_r^0$ and $\phi_r^1$
glue together to define a map of left $\O$-modules 
\[ \phi_r:\O(n-k)\rightarrow \Pr^k(\O (n)) .\]
For $\phi_r$ to be well-defined it is necessary that $\phi_r^0$ and 
$\phi_r^1$ agree
on $\U_{01}$. We have that $x_0^{n-k}=t^{k-n}x_1^{n-k}$, and this implies that
$\phi_r^0(x_0^{n-k})=\phi_r^0(t^{k-n}x_1^{n-k})=t^{k-n}\phi_r^1(x_1^{n-k})$ 
and we get from equation \ref{phi1def} 
\[ \phi_r^0(x_0^{n-k})=t^{k-n}\sum_{j=0}^rx_{j,r}t^{j-r}ds^j\otimes x_1^n\]
which equals
\begin{equation}\label{formula1} \sum_{j=0}^r x_{j,r}t^{j+k-n-r}ds^j\otimes x_1^n.
\end{equation}
Using relation \ref{basisrelation}, we substitute $ds^j\otimes x_1^{n}$ in 
formula \ref{formula1}
and get the following expression:
\[ \phi^0_r(x_0^{n-k})=\sum_{j=0}^r x_{j,r}t^{j+k-n-r}\left( \sum_{i=0}^{k-j}
(-1)^jt^{n-i-2j}\binom{n-j}{i}dt^{i+j}\otimes x_0^n\right).\]
Let  $l=i+j$ be a change of index. We get the expression
\[ \phi^0_r(x_0^{n-k})=\sum_{j=0}^r x_{j,r}t^{j+k-n-r}\left(\sum_{l=j}^{k}
(-1)^j t^{n-l-j}\binom{n-j}{l-j}dt^{l}\otimes x_0^n \right).\]
Since $\binom{a}{-1}=\binom{a}{-2}=\cdots =0$, we get
the expression
\[ \phi^0_r(x_0^{n-k})=\sum_{j=0}^r x_{j,r}t^{j+k-n-r}\left( \sum_{l=0}^{k}
(-1)^j t^{n-l-j}\binom{n-j}{l-j}dt^{l}\otimes x_0^n \right).\]
Simplify to obtain
\[ \sum_{j=0}^r \sum_{l=0}^k (-1)^j t^{k-r-l}\binom{n-j}{l-j}x_{j,r}
dt^l\otimes x_0^n .\]
Change order of summation and simplify:
\[ \phi_r^0(x_0^{n-k})=\sum_{l=0}^k t^{k-r-l} \left( \sum_{j=0}^r (-1)^j 
\binom{n-j}{l-j} x_{j,r} \right) dt^l\otimes x_0^n.\]
Let 
\begin{equation} \label{linearequations} c^r_l=\sum_{j=0}^r (-1)^j 
\binom{n-j}{l-j} x_{j,r}
\end{equation}
for $r=1,\dots, k$ and $l=0,\dots ,k$.  
We have defined maps 
\begin{equation} \label{def1}
 \phi^0_r(x_0^{n-k})=\sum_{l=0}^k t^{k-r-l}c^r_l dt^l\otimes x_0^n 
\end{equation}
and
\begin{equation} \label{def2}
 \phi^1_r(x_1^{n-k})=\sum_{j=0}^r x_{j,r}t^{j-r} ds^j\otimes x_1^n.
\end{equation}
Note that the definitions from \ref{def1} and \ref{def2} are 
valid for $r=0,\dots ,k$
since $c^0_i=\binom{n}{i}$ and $x_{0,0}=1$.
\begin{lemma} \label{maplemma}
Let $r=1,\dots, k$. The maps $\phi^0_r$ and $\phi^1_r$
glue to a well-defined map of left $\O$-modules
\[ \phi_r:\O(n-k)\rightarrow \Pr^1(\O(n)) \]
if and only if $c^r_k=c^r_{k-1}=\cdots =c^r_{k-r+1}=0$ and
$c^r_{k-r}=1$
\end{lemma}
\begin{proof} Consider the expression from \ref{def1}:
\[
\phi_r^0(x_0^{n-k})=\sum_{l=0}^{k-r-1}c^r_lt^{k-r-l}dt^l\otimes x_0^n+ 
c^{k-r}_r dt^{k-r}\otimes x_0^n +\]
\[c^r_{k-r+1}\t dt^{k-r+1}\otimes x_0^n+ c^r_{k-r+2}\frac{1}{t^2}dt^{k-r+2}\otimes 
x_0^n + \cdots + c^r_k\frac{1}{t^r}dt^k\otimes x_0^n.\]
We see that  the maps $\phi_r^0$ and $\phi_r^1$ glue if and only if
we have 
\begin{equation} \label{lineqs}
c^r_{k}=c^r_{k-1}=\cdots =c^r_{k-r+1}=0  \end{equation}
and 
\[ c^{r}_{k-r}=1,\]
and the lemma follows.
\end{proof}

We shall prove the following theorem: 
%
%
%
%
%
\begin{theorem} \label{maintheorem}
There exists $k$ systems $\{A_r\mathbf{x_r}=b_r\}_{r=1}^k$
of linear equations, where $A_r\mathbf{x_r}=b_r$ 
is a system of $k+1$ linear equations in $k+1$ variables, with coefficients in $\Z$. 
For all $r=1,\dots, k$, the system $A_r\mathbf{x_r}=b_r $ has the 
property that a solution $\mathbf{x_r}$  with coefficients in a field $F$, 
gives rise to 
an explicit map 
\[ \phi_r:\O(n-k)\rightarrow \Pr^k(\O (n)) \]
of left $\O$-modules.
Assume that there exists a field $F$ with the property that for all $r=1,\dots
,r$, 
$A_r\mathbf{x_r}=b_r$ has a solution $\mathbf{x_r}$ such that 
$\prod_{i=0}^kx_{i,i}\neq 0$, 
%
%
then the map
$\phi= \oplus_{i=0}^k \phi_i$ gives rise to an explicit isomorphism
\[ \Pr^k(\O (n))\cong \oplus_{i=0}^k \O(n-k) \]
of left $\O$-modules.
\end{theorem}
%
%
%
%
\begin{proof} Let $r=1,\dots , k$, and consider the equations from the proof
of Lemma \ref{maplemma}:
We have that $c_k^r=c_{k-1}^r=\cdots = c_{k-r+1}^r=0$ and that $c_{k-r}^r=1$.
We get from the equation $c_k^r=0$ the following:
\[\binom{n}{k}x_{0,r}-\binom{n-1}{k-1}x_{1,r}+\binom{n-2}{k-2}x_{2,r}+\cdots
+(-1)^r\binom{n-r}{k-r}x_{r,r}=0 .\]
Writing out $c_{k-1}^r=0$ we get the equation
\[\binom{n}{k-1}x_{0,r}-\binom{n-1}{k-2}x_{1,r}+\binom{n-2}{k-3}x_{2,r}+\cdots
+(-1)^r\binom{n-r}{k-r-1}x_{r,r}=0 .\]
The equation $c_{k-r}^r=1$ gives the following:
\[ \binom{n}{k-r}x_{0,r}-\binom{n-1}{k-r-1}x_{1,r}+\binom{n-2}{k-r-2}x_{2,r}+\cdots
+\binom{n-r}{k-2r}x_{r,r}=1.\]
We get a system of linear equations $A_r\mathbf{x_r}=b_r$, where $A_r$ 
is the 
rank $r+1$ matrix
\[
\begin{pmatrix} \binom{n}{k} & -\binom{n-1}{k-1} & \binom{n-2}{k-2} &  \cdots & (-1)^r\binom{n-r}{k-r} \\
\binom{n}{k-1} & -\binom{n-1}{k-1} & \binom{n-2}{k-3} &  \cdots & (-1)^r\binom{n-r}{k-r-1} \\
 \vdots & \vdots & \vdots & \vdots \\
\binom{n}{k-r+1} & -\binom{n-1}{k-r} & \binom{n-2}{k-r-1} &  \cdots 
& (-1)^r\binom{n-r}{k-2r+1} \\
\binom{n}{k-r} & -\binom{n-1}{k-r-1} & \binom{n-2}{k-r-2} &  \cdots & (-1)^r\binom{n-r}{k-2r} 
\end{pmatrix}
,\]
$\mathbf{x_r}$ is the rank $r+1$ vector $(x_{0,r},x_{1,r},\dots ,x_{r,r})$, and $b_r$
is the rank $r+1$-vector $(0,0,\dots ,0,1)$. Clearly the coefficients of $A_r$ and $b_r$ 
are in $\Z$. Also, assume $\mathbf{x_r}$ is a solution to the system $A_r\mathbf{x_r}=b_r$
with coefficients in a field $F$, then by construction and lemma 
\ref{maplemma}, the map 
\[ \phi_r:\O(n-k)\rightarrow \Pr^k(\O (n)) \] defined by
\[ \phi^0_r(x_0^{n-k})=\sum_{l=0}^{k-r}c^l_rt^{k-r-l}dt^l\otimes x_0^n \]
and 
\[ \phi^1_r(x_1^{n-k})=\sum_{j=0}^r x_{j,r}t^{j-r}ds^j\otimes x_1^n \]
is a well defined and nontrivial map of left $\O$-modules.
Assume that there exists a field $F$ and $k$ solutions 
$\mathbf{x_1},\dots ,\mathbf{x_k}$ to the 
systems $A_r\mathbf{x_r}=b_r$, with coefficients in $F$ satisfying the 
property that $\prod_{i=0}^k x_{i,i}\neq 0$. On the open set $\U_0$ the module
$\oplus_{i=0}^k\O(n-k)$ is a free  $k[t]$-module on the basis 
$\{x_0^{n-k}e_0,\dots ,x_0^{n-k}e_k\}$.
Define the map 
\[ \phi^0:\oplus_{i=0}^k \O(n-k)|_{\U_0}\rightarrow \Pr^k(\O(n))|_{\U_0} \]
as follows:
\[ \phi^0(x_0^{n-k}e_r)=\phi^0_r(x_0^{n-k}) .\]
On the open set $\U_1$ the module
$\oplus_{i=0}^k\O(n-k)$ is a free  $k[\t]$-module on the basis
$\{x_1^{n-k}f_0,\dots ,x_1^{n-k}f_k\}$. 
Define the map 
\[ \phi^1:\oplus_{i=0}^k \O(n-k)|_{\U_1}\rightarrow \Pr^k(\O(n))|_{\U_1} \]
as follows:
\[ \phi^1(x_1^{n-k}f_r)=\phi^1_r(x_1^{n-k}) .\] Then by construction, the
maps $\phi^0$ and $\phi^1$ glue to a well-defined map $\phi$ from
$\oplus_{i=0}^k \O(n-k)$ to $\Pr^k(\O(n))$ of left $\O$-modules. We
show explicitly that the map $\phi$ is an isomorphism:
Consider the matrix corresponding to the map $\phi|_{\U_0}=[\phi^0]$.
\[ [\phi^0]=
\begin{pmatrix} t^kc^0_0 & t^{k-1}c^1_0 & \cdots & tc^{k-1}_0 & c^k_0 \\
t^{k-1}c^0_1 & t^{k-2}c^1_1 & \cdots & c^{k-1}_1 & 0 \\
\vdots & \vdots & \cdots & \vdots & \vdots \\
tc^0_{k-1} & c^1_{k-1} & \cdots & 0 & 0 \\
c^0_k & 0 & \cdots & 0 & 0 \end{pmatrix} .\]
We see that the determinant $|[\phi^0]|$ equals $\prod_{i=0}^kc^i_{k-i}$
which equals $1$ by construction, hence the map $\phi^0$ is an isomorphism.
Consider the matrix $\phi^1|_{\U_1}=[\phi^1]$.
\[ [\phi^1]=
\begin{pmatrix} x_{0,0} & x_{0,1} \frac{1}{t} & \cdots & x_{0,k-1}
\frac{1}{t^{k-1}} & x_{0,k} \frac{1}{t^k} \\
0 & x_{1,1} & \cdots & x_{1,k-1}\frac{1}{t^{k-2}} & x_{1,k}\frac{1}{t^{k-1}} \\
\vdots & \vdots & \cdots & \vdots & \vdots \\
0 & 0 & \cdots & 0 & x_{k,k} \end{pmatrix} .\]
The determinant $|[\phi^1]|$ equals $\prod_{i=0}^kx_{i,i}$ which 
is non-zero by hypothesis, and the theorem follows.
\end{proof}

\section{Application: The left module-structure in characteristic zero}
In this section we use the results obtained in the previous section to 
determine
the splitting-type of $\Pr^k(\O(n))$ for all $1\leq k \leq n$ on the 
projective line defined over any field of characteristic zero.
\begin{lemma} \label{binomialreduction} 
Let $n,k,a,b\geq 0$, and put $\binom{0}{1}=0$. Then we have that
\[ \binom{n-a+1}{k-a-b+2}-\frac{n-k+b}{k-b+1}\binom{n-a+1}{k-a-b+1}=
\frac{ \binom{a-1}{1}\binom{n-a+1}{k-a-b+2} }{\binom{k-b+1}{1} }.\]
\end{lemma}
\begin{proof} Easy calculation.
\end{proof}
%
%
%
%
\begin{proposition} \label{determinants}
Let $n,k\geq 1$, and consider the matrix $A_r$ from Theorem \ref{maintheorem}
, with $r=1,2,\dots$. Then the determinant
\[ |A_r|=
\begin{vmatrix} \binom{n}{k} & -\binom{n-1}{k-1} & \binom{n-2}{k-2} &  \cdots & (-1)^r\binom{n-r}{k-r} \\
\binom{n}{k-1} & -\binom{n-1}{k-1} & \binom{n-2}{k-3} &  \cdots & (-1)^r\binom{n-r}{k-r-1} \\
 \vdots & \vdots & \vdots & \vdots \\
\binom{n}{k-r+1} & -\binom{n-1}{k-r} & \binom{n-2}{k-r-1} &  \cdots 
& (-1)^r\binom{n-r}{k-2r+1} \\
\binom{n}{k-r} & -\binom{n-1}{k-r-1} & \binom{n-2}{k-r-2} &  \cdots & (-1)^r\binom{n-r}{k-2r} 
\end{vmatrix}
\]
equals
\[ \prod_{l=0}^r\frac{\binom{n-l}{k-r} }
{\binom{k-l}{r-l}}\]
up to sign.
\end{proposition} 
%
%
\begin{proof} We prove the formula by induction on the rank of the matrix. 
Assume $r=1$. Add $-\frac{n-k+1}{k}$ times the second row to the first row 
of $A_1$, and apply Lemma \ref{binomialreduction} with $a=1$ and $b=1$, and we 
see that the formula is true for $r=1$. Assume the formula is true for rank
$r$, matrices $A_{r-1}$.  
Consider the matrix 
\[ M_r=
\begin{pmatrix} \binom{n}{k} & \binom{n-1}{k-1} & \binom{n-2}{k-2} &  \cdots & 
 \binom{n-r}{k-r} \\
\binom{n}{k-1} & \binom{n-1}{k-1} & \binom{n-2}{k-3} &  \cdots & 
\binom{n-r}{k-r-1} \\
 \vdots & \vdots & \vdots & \vdots \\
\binom{n}{k-r+1} & \binom{n-1}{k-r} & \binom{n-2}{k-r-1} &  \cdots 
& \binom{n-r}{k-2r+1} \\
\binom{n}{k-r} & \binom{n-1}{k-r-1} & \binom{n-2}{k-r-2} &  \cdots & 
\binom{n-r}{k-2r} 
\end{pmatrix}
\] 
which is the matrix $A_r$ with signs removed.
Add $-\frac{n-k+1}{k}$ times the second row to the first row. Continue and 
add $-\frac{n-k+1+i}{k-i}$ times the $i+1$'th row to the $i$'th row, for
$i=2,\dots, r-1$. Apply lemma \ref{binomialreduction} to get the matrix
\[ N_r=\begin{pmatrix}
0 & \frac{\binom{1}{1}\binom{n-1}{k-1} }{\binom{k}{1} } & 
\frac{ \binom{2}{1}\binom{n-2}{k-2}}{\binom{k}{1} } & \cdots &
\frac{ \binom{r}{1}\binom{n-r}{k-r} }{\binom{k}{1} } \\
0 & \frac{\binom{1}{1}\binom{n-1}{k-2} }{\binom{k-1}{1} } & 
\frac{ \binom{2}{1}\binom{n-2}{k-3}}{\binom{k-1}{1} } & \cdots &
\frac{ \binom{r}{1}\binom{n-r}{k-r-1} }{\binom{k-1}{1} } \\
0 & \frac{\binom{1}{1}\binom{n-1}{k-3} }{\binom{k-2}{1} } & 
\frac{ \binom{2}{1}\binom{n-2}{k-4}}{\binom{k-2}{1} } & \cdots &
\frac{ \binom{r}{1}\binom{n-r}{k-r-2} }{\binom{k-2}{1} } \\
\vdots & \vdots & \vdots & \vdots & \vdots \\
0 & \frac{\binom{1}{1}\binom{n-1}{k-r} }{\binom{k-r+1}{1} } & 
\frac{ \binom{2}{1}\binom{n-2}{k-r-1}}{\binom{k-r+1}{1} } & \cdots &
\frac{ \binom{r}{1}\binom{n-r}{k-2r+1} }{\binom{k-r+1}{1} } \\
\binom{n}{k-r} & * & * & \cdots & * 
\end{pmatrix}.\]
The determinant of $N_r$ equals 
\[ (-1)^r\frac{\binom{n}{k-r}\binom{1}{1}\binom{2}{1}\cdots \binom{r}{1} }
{\binom{k-r+1}{1}\binom{k-r+2}{1}\cdots \binom{k-1}{1}\binom{k}{1} }|M'_{r-1}|,\]
where $M'_{r-1}$ is the matrix
\[ \begin{pmatrix} \binom{n'}{k'} & \binom{n'-1}{k'-1} & \cdots &
\binom{n'-r+1}{k'-r+1} \\
\binom{n'}{k'-1} & \binom{n'-1}{k'-2} & \cdots &
\binom{n'-r+1}{k'-r} \\
\binom{n'}{k'-2} & \binom{n'-1}{k'-3} & \cdots &
\binom{n'-r+1}{k'-r11} \\
\vdots & \vdots & \cdots & \vdots \\
\binom{n'}{k'-r+1} & \binom{n'-1}{k'-r} & \cdots &
\binom{n'-r+1}{k'-2r+2} 
\end{pmatrix},\]
and $n'=n-1$ and $k'=k-1$.
By the induction hypothesis, we get modulo signs
\[ \frac{\binom{n}{k-r}}{\binom{k}{r}} \prod_{i=0}^{r-1}
\frac{\binom{n'-i}{k'-r+1}}{\binom{k'-i}{r-1-i}}.\]
Change index by letting $l=i+1$ we get
\[ \frac{\binom{n}{k-r}}{\binom{k}{r}} \prod_{i=1}^{r}
\frac{\binom{n-l}{k-r}}{\binom{k-l}{r-l}}=
\prod_{l=0}^{r}
\frac{\binom{n-l}{k-r}}{\binom{k-l}{r-l}}\]
and the proposition follows.
\end{proof}
%
%
%
%
\begin{proposition} \label{characteristiczero}
Let $F$ be a field of characteristic zero. Then the maps
$\phi_r$ from Theorem \ref{maintheorem} exist for $r=1,\dots, k$, 
and the induced map 
$\phi=\oplus _{i=0}^r \phi_r$ defines an isomorphism
\[ \Pr^k(\O (n))\cong \oplus_{i=0}^k\O(n-k) \]
of left $\O$-modules.
\end{proposition}
%
%
\begin{proof} Consider the systems $A_r\mathbf{x_r}=b_r$ for $r=1,\dots,k$, from 
the proof of theorem \ref{maintheorem}. By Proposition \ref{determinants}
we have that 
\[ |A_r|= \prod_{l=0}^r\frac{ \binom{n-l}{k-r} }{
\binom{k-l}{r-l} } \]
modulo signs.
Since the characteristic of $F$ is zero, it follows that $|A_r|\neq 0$
for all $r=1,\dots ,k$, hence the system $A_r\mathbf{x_r}=b_r$ has a unique 
solution $\mathbf{x_r}=A_r^{-1}b_r$ for all $r$. It follows from 
theorem \ref{maintheorem} that maps
\[ \phi_r:\O(n-k)\rightarrow \Pr^k(\O(n)) \]
of left $\O$-modules exist, for $r=1,\dots ,k$, and we can consider the map
\[ \phi=\oplus_{i=0}^r\phi_i:\oplus_{i=0}^r\O(n-k)\rightarrow \Pr^k(\O(n)) .\]
We want to prove that $\phi$ is an isomorphism. Again by Theorem \ref{maintheorem}
$\phi$ is an isomorphism if and only if $x_{i,i}\neq 0$ for $i=0,\dots ,k$.
Assume that $x_{r,r}=0$, and consider 
the system $A_r\mathbf{x_r}=b_r$. If $x_{r,r}=0$ we get a new system
$A_{r-1}\mathbf{y_{r-1}}=0$, where $\mathbf{y_{r-1}}$ is the vector
$(x_{0,r},\dots ,x_{r-1,r})$. Since the matrix $A_{r-1}$ is invertible, 
it follows that the system $A_{r-1}\mathbf{y_{r-1}}=0$ only has the trivial
solution $\mathbf{y_{r-1}}=0$, hence $x_{0,r}=\cdots =x_{r-1,r}=0$, and we have arrived at a contradiction to the assumption that $\mathbf{x_r}$ 
is a solution to the system $A_r\mathbf{x_r}=b_r$, where $b_r$ is the vector
$(0,\dots,0,1)$. It follows that $x_{r,r}\neq 0$ for all $r=0,\dots ,k$, 
and the proposition follows from Theorem \ref{maintheorem}.
\end{proof}

\section{Splitting the right module-structure}
In this section we consider the splitting type of the principal parts as
left and right $\O$-module on $\P^1$ defined over $F$, where $F$ is any field.
We prove that in most cases the splitting-type
as left module differs from the splitting-type as right module.
We also show how the splitting-type of the principal parts as
left $\O_{\P^1}$-module changes with the characteristic of the field
$F$.
%
%
%
\begin{theorem} \label{leftandright} 
Consider $\Pr^1(\O(n))$ on $\P^1_F$ where $F$ is any field
and $n\geq 1$. Then $\Pr^1(\O(n))$ is locally free as left and right $\O$-module.
If the characteristic of $F$ does not divide $n$, $\Pr^1(\O(n))$ splits as 
$\O(n-1)\oplus \O(n-1)$ as left $\O$-module and as $\O(n)\oplus \O(n-2)$ as right
$\O$-module. If the characteristic of $F$ divides $n$, $\Pr^1(\O(n))$ splits
as $\O(n)\oplus \O(n-2)$ as left and right $\O$-module.
\end{theorem}
\begin{proof} Recall from Section $4$, that the splitting type of
$\Pr^1(\O (n))$ as left $\O$-module is $\O(n-1)\oplus \O(n-1)$ if the characteristic of $F$ does not divide $n$. We next consider the right $\O$-module structure. Let $p$ and $q$ be the projection maps from $\P^1 \times \P^1$ to $\P^1$.
By definition, $\Pr^1(\O(n))$ is $p_*(\O_{\Delta^1}\otimes q^*\O(n))$, where
$\O_{\Delta^1}$ is the first order infinitesimal neighborhood of the 
diagonal. By Proposition \ref{pushiso}
we get the right $\O$-module structure, by considering the module
$q_*(\O_{\Delta^1}\otimes q^*\O(n))=q_*(\O_{\Delta^1})\otimes \O(n)$. And we 
see that $\E=q_*(\O_{\Delta^1})\otimes \O(n)$ is locally free, since it is 
the tensor-product of two locally free sheaves. One checks that 
$\E|_{\U_0}$ is a free
$k[u]$-module on the basis $E=\{1\otimes x_0^n, du\otimes x_0^n\}$ where
$du=t-u$. Similarly,
$\E|_{\U_1}$ is a free $k[\u]$-module on the basis $E'=\{1\otimes x_1^n,
d(\u)\otimes x_1^n\}$ where $d(\u)=\t-\u$. On $\U_{01}$ the module $\E$ is free on $E$ and $E'$ as 
$F[u,\u]$-module. We compute the transition-matrix $[R]^{E'}_{E}$.
We see that $1\otimes x_1^n$ equals $u^n\otimes x_0^n$. 
By definition, $d(\u)\otimes x_1^n$ equals $(\t-\u)u^n\otimes x_0^n$.
We get
\[ (\t-\u)u^n\otimes x_0^n=\frac{u-t}{ut}u^n\otimes x_0^n=
-u^{n-1}du(\t)\otimes x_0^n.\]
By definition, $t=u+du$, hence we get
\[ -u^{n-1}du\frac{1}{u+du}\otimes x_0^n=-u^{n-2}du\frac{1}{1+du/u}\otimes x_0^n.\]
And since $du^2=du^3=\cdots =0$, we see that
\[d(\u)\otimes x_1^n=-u^{n-2}du\otimes x_0^n,\]
hence the transition-matrix looks as follows:
\[ [R]^{E'}_E=
\begin{pmatrix}  u^n & 0 \\
      0 & -u^{n-2} \end{pmatrix} ,\]
and it follows that $\E$ splits as $\O(n)\oplus \O(n-2)$ as $\O$-module.
Recall the transition-matrix for $\Pr^1(\O(n))$ as left $\O$-module,
\[ [L]^{C'}_C=
\begin{pmatrix} t^n & 0 \\ 
 nt^{n-1} & -t^{n-2} \end{pmatrix}. \]
Then we see that if the characteristic of $F$ divides $n$, the splitting-type
of $\Pr^1(\O(n))$ is $\O(n)\oplus \O(n-2)$ as left and right $\O$-module, 
and we have proved the theorem.
\end{proof}

\end{document}